\documentclass[12pt]{article}
\usepackage{amssymb}
\usepackage{amsmath}
\usepackage{subfig}
\usepackage{pgfplots}
\usepackage{natbib}
\newenvironment{proof}{\noindent{\bf Proof: }}{\hfill\fbox{}\vspace*{3mm}}

\def\b{\bf}

\def\pa{\partial}

\def\dfrac{\displaystyle \frac}
\def\dsum{\displaystyle \sum}

\def\dprod{\displaystyle \prod}
\def\b1{{\bf{1}}}
\def\={\!\!\!=\!\!\!}

\begin{document}
\title{\textbf{An EOQ model for imperfect quality items with multiple screening and shortage backordering}}

\author{Allen H. Tai\\
Department of Applied Mathematics,\\
The Hong Kong Polytechnic University,\\
Hung Hom, Hong Kong}
\date{}
\maketitle
\begin{abstract}
In this paper, we propose an inventory model where items are inspected through multiple screening processes before delivery to customers.
Each screening process has independent screening rate and defective percentage.
Defective items screened out are stored and then returned to supplier. 
Shortage backordering are also allowed in the model.
Two approaches are used to obtain the closed-form optimal order size and the maximum backordering quantity.
Numerical examples are also provided to demonstrate the use of the model.
\end{abstract}

\section{Introduction}
When a batch of items arrives from supplier, imperfect quality items do exist in the batch.
Therefore, screening processes are required before the items are delivered to customers.
\citet{Salameh} first extended the classical economic order quantity (EOQ) model by integrating a screening process for imperfect quality items into the model.
Since then many researchers proposed extensions to the model developed by \citet{Salameh}.
\citet{Goyal} proposed an approximation to the model by simplifying the expected revenue and expected cost per unit time.
\citet{Maddah} pointed out a flaw in \citet{Salameh} and reformulated the model using the renewal reward theorem.
\citet{Wahab} considered the case when there are different holding costs for the imperfect quality and screened items.
Other extensions of this original model can be found in the review paper by \citet{Khan}.

One natural extension to an EOQ model is to consider shortage situation.
\citet{Wee} developed an inventory model for items with imperfect quality and shortage backordering.
\citet{Chang} then revisited their study and provided closed-form solutions using the renewal reward theorem.
However, in both works they assumed that the items satisfying the backorder are delivered without any screening process.
\citet{Hsu} then proposed a model which items are first screened before delivery to customers.
In this paper, we aim at providing a general version of the work of \citet{Hsu} by considering multiple screening processes and shortage backordering.

Suppose that after the replenishment of items arrives from the supplier, the items have to go through $n$ screening tests before delivering to customers. 
Denote by $S_i$ the type $i$ screening process and by $x_i$ the corresponding screening rate of $S_i$.
Without loss of generality, we suppose
$$
x_1 \ge x_2 \ge \cdots \ge x_n.
$$
 
Let $p_i$ be the proportion of items in the lot that can not pass $S_i$. 
Here we assume that all $p_i$ are independent of each others in the sense that whether an item can pass $S_i$ does not depend on the result of other screening processes.
Since for any $i<j$, $S_i$ will finish before $S_j$, 
the proportion of items in $y$ that is screened out after $S_i$ is 
$$\rho_i=
\left\{
\begin{array}{ll}
p_1, &\mbox{for }i=1;\\[2mm]
\dprod_{k=1}^{i-1} (1-p_k)p_i, &\mbox{for }i=2,\ldots,n.
\end{array}
\right.
$$
Then the total proportion of items in $y$ that are classified as defective is given by
$$
\rho=\dsum_{i=1}^n \rho_i.
$$
The items go through screening processes one by one.
Here we assume that the defective items are evenly distributed in the lot.
Hence, the screening time for the $S_i$ can be modeled as $y/x_i$.
The screening processes are illustrated in Figure 1.
In Figure 1, the arrows represent the screening rates and the lines represent the items undergoing the corresponding screening process. 

The remainder of the paper is organized as follows. 
In Section 2, we present the mathematical model for the problem and derive the expected profit per unit time using the renewal reward theorem. 
In Section 3, we obtain the optimal order size and the maximum backordering quantity with an analysis of the model.
We then give numerical examples to demonstrate the use of our
model in Section 4. 
Finally, concluding remarks are given in Section 5 to conclude the paper
and address further research issues.

\section{Mathematical model}

In this section, we first give the notations and assumptions used in the development of the model.
Some of the notations and assumptions were also used in \citet{Hsu} and \citet{Wee}.

\noindent Notations:

\begin{tabular}{ll}
$y$ & the order size,\\
$D$ & the demand rate,\\
$n$ & number of screening processes in a cycle,\\
$S_i$ & type $i$ screening process ($i=1,\ldots,n$),\\
$x_i$ & the screening rate of $S_i$,\\
$c$ & the purchasing cost per unit,\\
$K$ & the ordering cost per order,\\
$p_i$ & the defective percentage of $S_i$ in $y$\\
$\rho_i$ & the proportion of items in $y$ that screen out after type $i$
screening process \\
&($i=1,\ldots,n$),\\
$s$ & the selling price of good quality items, \\
$v$ & the salvage value of defective items ($v<s$),\\
$d_i$ & the unit screening cost of $S_i$,\\
$B$ & the maximum backordering quantity,\\
$b$ & the backordering cost per unit per unit time,\\
$h$ & the holding cost per unit per unit time,\\
$h_d$ & the holding cost for defective items per unit per unit time,\\
$^*$ & the superscript representing optimal value.

\end{tabular}

\noindent Assumptions:
\begin{enumerate}
\item The demand rate is a known constant.
\item The lead time of inventory replenishment is assumed to be negligible.
\item All screening processes and demands proceed simultaneously, 
but the screening rates are greater than the demand rate,
$x_i >D$ ($i=1,\ldots,n$).
\item The defective items exist in lot size $y$.
	\item To avoid shortages within screening time, we assume that
	$$
	\rho \leq 1-\dfrac{D}{x_n}.
	$$
	\item The defective items are returned to the supplier when the replenishment items arrive.
	\item Shortage is completely backordered.
	\item A single product is considered.
\end{enumerate}

The inventory level in a replenishment cycle is illustrated in Figure 2.
At the beginning of the replenishment cycle, all screening processes proceed simultaneously.
Therefore, the rate of items to complete all the screening processes depends only on the lowest screening process rate $x_n$.
After completing all the screening processes, the items were first shipped to satisfy the demands. 
Hence, the inventory level decreases at a rate of $(1-\rho)x_n$.
The remaining items are then used to clear the outstanding backorders.
Hence, the backordering quantity decreases at a rate of $(1-\rho)x_n-D$.
The screening time for $S_i$ is $y/x_i$. 
After the screening process $S_i$, the defective items ($\rho_i y$ units) are transferred to another inventory warehouse.
The inventory level of defective items is illustrated in Figure 3.
The aim of this paper is to develop a mathematical model such that the optimal order size and the maximum backordering quantity can be obtained.
We first give the length of time intervals in our model.
\begin{enumerate}
	\item[(i)] The replenishment cycle length $T$: $T=\dfrac{(1-\rho)y}{D}$;
	\item[(ii)] $t_1=\dfrac{B}{D}$;
  \item[(iii)] $t_2=T-t_1=\dfrac{(1-\rho)y-B}{D}$;
  \item[(iv)] $t_3=\dfrac{B}{(1-\rho)x_n-D}$;
	\item[(v)] $t_4=t_2-t_3=\dfrac{(1-\rho)y-B}{D}-\dfrac{B}{(1-\rho)x_n-D}$.
\end{enumerate}
Denote by $TR(y)$ and $TC(y,B)$ the total revenue and the total cost per cycle respectively. 
$TR(y)$ is the sum of the total sales of good quality items and the amount received from the supplier for the return of the imperfect quality items. 
One has
$$TR(y) = (1 - \rho)ys + \rho yv.$$
$TC(y,B)$ consists of the following costs:
\begin{enumerate}
\item[(i)] The purchasing cost: $cy$;
\item[(ii)] The screening cost:
		$$
		 \dsum_{i=1}^n \Big[d_i\times\Big(1-\dsum_{k=0}^{i-1} \rho_k\Big)y\Big]=\Big[\dsum_{i=1}^n d_i-\dsum_{i=1}^{n-1}\rho_i  \Big( \dsum_{k=i+1}^n d_k \Big)\Big]y;
		$$

\item[(iii)] The holding cost for good quality items:
	$$
	\begin{array}{cl}
	&h\times  \Big[ \dfrac{t_4^2\times D}{2}+\dfrac{t_3^2\times (1-\rho)x_n}{2}
	+\dsum_{i=1}^n \rho_i y\times \dfrac{y}{x_i}+t_3\times t_4\times D\Big]\\
	=&h  \Big[ \dfrac{(t_2-t_3)^2 D}{2}+\dfrac{t_3^2 (1-\rho)x_n}{2}
	+\dsum_{i=1}^n \rho_i y \dfrac{y}{x_i}+t_3 (t_2-t_3) D\Big]\\
	=&h  \Big[ \dfrac{t_2^2 D}{2}+\dfrac{t_3^2 \big((1-\rho)x_n-D\big)}{2}
	+\dsum_{i=1}^n \rho_i y \dfrac{y}{x_i}\Big]\\
	=&h  \Big[ \dfrac{\big((1-\rho)y-B\big)^2}{2D}+\dfrac{B^2}{2\big((1-\rho)x_n-D\big)}
	+\dsum_{i=1}^n \rho_i y \dfrac{y}{x_i}\Big]\\

	\end{array}
	$$
	
\item[(iv)] The holding cost for defective items:
		$$
	\begin{array}{rcl}
	h_d \times \Big[ \dsum_{i=1}^n \rho_i y \times \Big( T- \dfrac{y}{x_i}\Big)\Big]
	&=& h_d \Big[ \dsum_{i=1}^n \rho_i \Big( \dfrac{(1-\rho)}{D}- \dfrac{1}{x_i}\Big)
		\Big]y^2\\
		&=& h_d \Big[  \Big( \dfrac{\rho(1-\rho)}{D}- \dsum_{i=1}^n\dfrac{\rho_i}{x_i}\Big)
		\Big]y^2;
	\end{array}
	$$
	
\item[(v)] The shortage cost:
	$$
		\begin{array}{rcl}
	b \times \Big( \dfrac{t_3\times B}{2}+ \dfrac{t_1\times B}{2} \Big) 
	&=&	b \Big[ \dfrac{1}{2\big((1-\rho)x_n-D\big)}+ \dfrac{1}{2D} \Big]B^2\\[5mm]
		&=&	b \Big[ \dfrac{1-\rho}{2D\big((1-\rho)-D/x_n\big)} \Big]B^2;
		\end{array}	
		$$
		\item[(vi)] The ordering cost: $K$.
\end{enumerate}

For notational convenience, let 
$$P_1=E[(1-\rho)^2],\ P_2=E[\rho(1-\rho)],\ P_3=\dfrac{1}{1-E[\rho]},\ P_4=\dfrac{E[\rho]}{1-E[\rho]},$$
which depend on $\rho$ only. Let 
$$R=\dfrac{1-E[\rho]}{E\Big[ \dfrac{1-\rho}{(1-\rho)-D/x_n} \Big]},$$
which depends on $\rho, D$ and $x_n$. Let 
$$A_1=\dsum_{i=1}^n \dfrac{E[\rho_i]}{x_i},$$ 
which depends on $\rho_i, x_i$ $(i=1,\ldots,n)$ and  
$$A_2=\dsum_{i=1}^n d_i-\dsum_{i=1}^{n-1}E[\rho_i]  \Big( \dsum_{k=i+1}^n d_k \Big),$$ 
which depends on $\rho_i, d_i$ $(i=1,\ldots,n)$.
Then the expected net profit per cycle is given by
$$
\begin{array}{cl}
&ETP(y,B)\\
=&(1 - E[\rho])ys + E[\rho] yv-\Bigg\{cy+A_2y\\
&+h\Big[\dfrac{P_1}{2D}y^2-\dfrac{1 - E[\rho]}{D}By+\dfrac{B^2}{2RP_3D}  +A_1 y^2\Big]
+h_d \Big[  \Big( \dfrac{P_2}{D}- A_1\Big) \Big]y^2\\
&+\dfrac{b}{2RP_3D} B^2 +K \Bigg\}
\end{array}
$$
and the expected replenishment cycle length $E[T]=\dfrac{(1-E[\rho])y}{D}$.

By the renewal reward theorem, the expected profit per unit time is
\begin{equation}\label{ETPU}
\begin{array}{rcl}
ETPU(y,B)&=&\dfrac{ETP(y,B)}{E[T]}\\
&=&sD + vDP_4-(c+A_2)DP_3-\Bigg\{\dfrac{(h+b)B^2}{2Ry}-h B\\
&&+\Big[\dfrac{hP_1}{2}+hA_1 D+h_d(P_2-A_1 D) \Big]P_3y+\dfrac{KDP_3}{y}\Bigg\}.
\end{array}
\end{equation}
Our aim is to find the optimal order size $y$ and the maximum backordering quality $B$ such that the expected profit per unit time $ETPU(y,B)$ is maximized.
Since the first three terms in (\ref{ETPU}) are independent of $B$ and $y$,
the optimization problem reduces to minimizing $f(y,B)$ which is given by
\begin{equation}\label{fBy}
\begin{array}{rcl}
f(y,B)&=& \dfrac{(h+b) B^2}{2Ry}-hB \\[2mm]
&&+\Big[\dfrac{hP_1}{2}+hA_1 D+h_d(P_2-A_1 D) \Big]P_3y+\dfrac{KDP_3}{y}.
\end{array}
\end{equation}

\section{Analysis}
In this section, we present two different approaches to minimize $f(y,B)$ in (\ref{fBy}).

\noindent \textit{Approach 1}

We use an approach similar to \citet{Chang} here.
Notice that (\ref{fBy}) can be rewritten as
\begin{equation}\label{fBy1}
\begin{array}{rcl}
f(y,B)&=& \dfrac{h+b}{2Ry}\Big(B-\dfrac{hR}{h+b}y\Big)^2\\
&&+\Big[hP_1+2hA_1 D+2h_d(P_2-A_1 D)-\dfrac{h^2 R}{(h+b)P_3} \Big]\dfrac{P_3y}{2}\\[2mm]
&&+\dfrac{KDP_3}{y}.
\end{array}
\end{equation}
Hence, for any fixed $y$, $f(y,B)$ is minimized when
\begin{equation}\label{Biny}
B=\dfrac{hR}{h+b}y.
\end{equation}
Hence the problem reduces to minimizing
\begin{equation}\label{fy}
f(y)= \Big[ hP_1+2hA_1 D+2h_d(P_2-A_1 D)-\dfrac{h^2 R}{(h+b)P_3} \Big] \dfrac{P_3y}{2}+\dfrac{KDP_3}{y}.
\end{equation}
Applying the AM-GM inequality gives
$$
f(y) \geq P_3\sqrt{2KD \Big[hP_1+2hA_1 D+2h_d(P_2-A_1 D)-\dfrac{h^2 R}{(h+b)P_3} \Big]}.
$$
The equality holds when
$$\Big[hP_1+2hA_1 D+2h_d(P_2-A_1 D)-\dfrac{h^2 R}{(h+b)P_3} \Big] \dfrac{P_3y}{2}=\dfrac{KDP_3}{y},$$
which means $f(y)$ is minimized when
\begin{equation}\label{ystar}
y^*=\sqrt{\dfrac{2KD}{hP_1+2hA_1 D+2h_d(P_2-A_1 D)-\dfrac{h^2 R}{(h+b)P_3}}}.
\end{equation}
From (\ref{Biny}), we have
\begin{equation}\label{Bstar}
B^*=\dfrac{h R}{h+b}y^*.
\end{equation}

Hence the minimum value of $f(y,B)$ is
$$
f(y^*,B^*)=P_3\sqrt{2KD \Big[hP_1+2hA_1 D+2h_d(P_2-A_1 D)-\dfrac{h^2 R}{(h+b)P_3} \Big]}
$$
and the maximum expected profit per unit time is
\begin{equation}\label{ETPU*}
\begin{array}{rcl}
ETPU^*&=&sD + vDP_4-(c+A_2)DP_3\\[2mm]
&&-P_3\sqrt{2KD \Big[hP_1+2hA_1 D+2h_d(P_2-A_1 D)-\dfrac{h^2 R}{(h+b)P_3} \Big]}.
\end{array}
\end{equation}

\noindent \textit{Approach 2}

Differentiating $f(y,B)$ in (\ref{fBy}) with respect to $B$ and $y$ respectively give
$$
\dfrac{\pa f}{\pa B} = \dfrac{(h+b)B}{Ry}-h,
$$
$$
\dfrac{\pa f}{\pa y} = 
 -\dfrac{(h+b) B^2}{2 R y^2}+\Big[\dfrac{hP_1}{2}+hA_1 D+h_d(P_2-A_1 D) \Big]P_3-\dfrac{KDP_3}{2y^2}.
$$
The second order partial derivatives are
$$
\dfrac{\pa^2 f}{\pa B^2} = \dfrac{h+b}{Ry},
$$
$$
\dfrac{\pa^2 f}{\pa y^2} = \dfrac{(h+b) B^2}{R y^3}+\dfrac{2KDP_3}{y^3},
$$
$$
\dfrac{\pa^2 f}{\pa B\pa y} = -\dfrac{(h+b) B}{Ry^2}.
$$
Notice that 
$$
\dfrac{\pa^2 f}{\pa B^2}\dfrac{\pa^2 f}{\pa y^2}-\Big( \dfrac{\pa^2 f}{\pa B\pa y} \Big)^2=
\dfrac{2(h+b)KDP_3}{Ry^4}.
$$
Since we assume that $\rho<1-D/x_n$, we have 
$\dfrac{\pa^2 f}{\pa B^2} >0 $ and 
$\dfrac{\pa^2 f}{\pa B^2}\dfrac{\pa^2 f}{\pa y^2}-\Big( \dfrac{\pa^2 f}{\pa B\pa y} \Big)^2>0$.
This implies that $f(y,B)$ is strictly convex for positive $B$ and $y$.
Hence the unique global minimum for positive $B$ and $y$ can be obtained by solving $\dfrac{\pa f}{\pa B} = 0$ and
$\dfrac{\pa f}{\pa y} = 0$, which gives the same results as in (\ref{ystar}) and (\ref{Bstar}).

Based on the results above, we then present some special case analyses.
\begin{enumerate}
	\item[(i)] From (\ref{ETPU*}), the defective items are with salvage value $vD P_4$ and the holding cost for the defective items is
	$2P_3\sqrt{KDh_d(P_2-A_1D)}$. Hence there will be no benefit to store the defective items if
	$$
	v \leq \dfrac{2}{E[\rho]}\sqrt{Kh_d\Big(\dfrac{P_2}{D}-A_1\Big)}. 
	$$
	\item[(ii)] If the supplier does not accept returns, i.e.\ $v=0$, then there is no storage for the defective items ($h_d=0$). 
	The unit selling price should be higher than
	$$
	s_{min}=P_3\Big[ c+A_2+\sqrt{\dfrac{2Kh}{D}\Big( P_1+2A_1D-\dfrac{hA_1}{h+b} \Big)} \Big],
	$$
	so that positive profit can be guaranteed, i.e.\ $ETPU^*>0$.
	
\end{enumerate}

We remark that in the literatures mentioned below, the defective items are not stored for return.
Hence we consider $h_d=0$ in the following cases for comparison.

\begin{enumerate}
	\item[(iii)] When $n=1$, the proposed model reduces to the one presented in \citet{Hsu},
	with 
	\begin{equation}\label{ystar1}
	y^*=\sqrt{\dfrac{2KD}{h\Big[P_1+2A_1 D-\dfrac{h R}{(h+b)P_3}\Big]}}
	\end{equation}
	and
	\begin{equation}\label{Bstar1}
	B^*=\dfrac{h R}{h+b}y^*.
	\end{equation}
	We remark that the above expressions are simpler than those obtained in \citet{Hsu}.
	As stated in \citet{Hsu},
	if the defective percentage $p_1$ follows uniform distribution with probability density function
	$$
	f_{p_1}(p_1)=\left\{
	\begin{array}{ll}
	\dfrac{1}{\beta}, \quad &\mbox{for } 0< p_1 < \beta;\\
	0, & \mbox{otherwise}.
	\end{array}\right.
	$$
	then 
	$$
	E[p_1]=\dfrac{\beta}{2},\quad E[(1-p_1)^2]=1-\beta+\dfrac{\beta^2}{3},
	$$
	$$
	E\Big[\dfrac{1-p_1}{1-p_1-D/x_1}\Big]=1+\dfrac{D}{\beta x_1}\ln\Big(\dfrac{1-D/x_1}{1-\beta-D/x_1}\Big).
	$$
	In what follows, we give an analysis which is absent in \citet{Hsu}.
	The equations (\ref{ystar1}) and (\ref{Bstar1}) can be reduced to
  \begin{equation}\label{ystar2}
	y^*=\sqrt{\dfrac{2KD}{h\Bigg\{1-\beta+\dfrac{\beta^2}{3}+\dfrac{D\beta}{x_1}-\dfrac{h (1-\beta/2)^2}{(h+b)}\Big[1+\dfrac{D}{\beta x_1}\ln\Big(\dfrac{1-D/x_1}{1-\beta-D/x_1}\Big)\Big]^{-1}\Bigg\}}}
	\end{equation}
	and
	\begin{equation}\label{Bstar2}
	B^*=\dfrac{h (2-\beta)}{2(h+b)}\Big[1+\dfrac{D}{\beta x_1}\ln\Big(\dfrac{1-D/x_1}{1-\beta-D/x_1}\Big)\Big]^{-1}y^*.
	\end{equation}
	For the case when $\beta$ is small, by the Taylor series expansion, we can use the approximation
	$$
	\ln\Big(\dfrac{1-D/x_1}{1-\beta-D/x_1}\Big)\approx \dfrac{\beta}{1-\beta-D/x_1}.
	$$
	By neglecting the terms with $\beta^2$, we can further reduce (\ref{ystar1}) and (\ref{Bstar1}) to
  \begin{equation}\label{y1star}
	y_1^*=\sqrt{\dfrac{2KD(h+b)}{h[h(1-\beta)D/x_1+b(1-\beta-D\beta/x_1)]}}
	\end{equation}
	and
  \begin{equation}\label{B1star}
		B_1^*=\dfrac{h (2-\beta)(1-\beta-D/x_1)}{2(h+b)(1-\beta)} y_1^*.
	\end{equation}
	
	Finally, if all products are of good quality, which means the screening rate can be set as $x_1\to \infty$,
	then we have
	$E[(1-\rho)^2]=1$, $E[\rho]=0$ and $E\Big[\dfrac{1-\rho}{1-\rho-D/x_1}\Big]=1$.
	Hence $y^*$ reduces to the classical EOQ with shortages.
	
	\item[(iv)] When $b \to \infty$, i.e.\ no shortages are allowed, then
	$$
	y^* \to \sqrt{\dfrac{2KD}{h(P_1+2 A_1 D)}},
	$$
	which is a generalization of the result obtained in \citet{Maddah}.
\end{enumerate}

\section{Numerical examples}
In this section, we present some numerical examples to demonstrate the use of our model.
We apply the following parameters which are also used in \citet{Hsu} and \citet{Wee}.

\begin{tabular}{ll}
Demand rate, $D$ &= 50,000 units/year,\\
Ordering cost, $K$ &= 100/cycle,\\
Holding cost, $h$ &= \$5/unit/year,\\
Purchasing cost, $c$ &= \$25/unit,\\
Backordering cost, $b$ &= \$10/unit/year,\\
Selling price of good quality items, $s$ &= \$50/unit,\\
The salvage value of defective items, $v$ &= \$20/unit.\\
\end{tabular}

Suppose that there are seven types of screening processes, $S_i$ $(i=1,\ldots,7)$.
The corresponding defective percentage, $p_i$ $(i=1,\ldots,7)$, are assumed to be uniformly distributed with respective pdf
	$$
	f_{p_i}(p_i)=\left\{
	\begin{array}{ll}
	\dfrac{1}{\beta_i}, \quad &\mbox{for } 0< p_i < \beta_i;\\
	0, & \mbox{otherwise}.
	\end{array}\right.
	$$
The parameter $\beta_i$, the screening rate $x_i$ and the screening cost $d_i$ for the screening process $S_i$ $(i=1,\ldots,7)$ are given in Table 1.
\begin{table}
\centering
\begin{tabular}{llllllll}
\hline
&$S_1$&$S_2$&$S_3$&$S_4$&$S_5$&$S_6$&$S_7$ \\
\hline
$\beta_i$ &0.01 &0.04 &0.1 & 0.01 & 0.04 &0.04 &0.1 \\
$x_i$ (unit/min) &1&1&1&2&2&0.5&0.5\\
$d_i$ (\$/unit) &0.5 &0.5 &0.5&1 &0.5 &1&0.3\\
\hline
\end{tabular}
\caption{Parameters for the screening processes $S_i$ $(i=1,\ldots,7)$.}
\end{table}
The optimal $y^*$ and $B^*$ from (\ref{ystar}) and (\ref{Bstar}) and the corresponding expected profit per unit time for each of the screening process are given in Table 2.
We also give the approximated $y_1^*$ and $B_1^*$ from (\ref{y1star}) and (\ref{B1star}) to demonstrate the effectiveness of the approximations.
The results show that the approximations are good as the errors of the expected profit per unit time are all less than $0.01\%$. 
The results of $S_1,S_2$ and $S_3$ show that when $\beta$ increases, $y^*$ increases but $B^*$ and $ETPU^*$ decreases.
The results of $S_2$ and $S_5$ show that when $x$ increases, all $y^*$, $B^*$ and $ETPU^*$ increases.
The results of $S_3$ and $S_7$ show that when both $\beta$ and $d$ increase, both $y^*$ and $B^*$ increase but $ETPU^*$ decreases.

\begin{table}
\centering
\scriptsize
\begin{tabular}{lllllllll}
\hline
&$S_1$&$S_2$&$S_3$&$S_4$&$S_5$&$S_6$&$S_7$\\
\hline
$y^*$	&	1624.85	&	1638.40	&	1664.90	&	1679.81	&	1699.16	&	1534.16	&	1542.35	\\
$B^*$	&	384.34	&	379.32	&	368.63	&	477.24	&	474.22	&	209.17	&	194.28	\\
$ETPU(y^*,B^*)$	&	1217432.76	&	1213159.67	&	1204203.81	&	1192509.48	&	1213382.37	&	1187226.30	&	1214227.78	\\
$y_1^*$	&	1630.52	&	1662.38	&	1732.15	&	1682.95	&	1712.65	&	1573.85	&	1651.62	\\
$B_1^*$	&	384.90	&	381.61	&	374.58	&	477.72	&	476.31	&	208.45	&	191.32	\\
$ETPU(y_1^*,B_1^*)$	&	1217432.72	&	1213158.94	&	1204198.33	&	1192509.47	&	1213382.17	&	1187223.69	&	1214208.42	\\
\hline
\end{tabular}
\caption{Numerical results for the screening processes $S_i$ $(i=1,\ldots,7)$.}
\end{table}
\normalsize

We then provide the results for different combinations of the screening processes. 
First we consider the case when there are two screening processes in the cycle.
The results are given in Table 3.
We then compare these results with those in Table 2.
We remark that for the case of two screening processes,
the optimal order size $y^*$ is between those in the cases of one screening process.
However, when there are two screening processes in the cycle the maximum backordering quantity $B^*$
is less than both $B^*$ in the cases of one screening process.

Next we consider the case when there are three screening processes in the cycle.
The results are given in Table 4.
We observe that, for example in the case of $S_1+S_4+S_6$,
the optimal order size $y^*$ is between those in the cases of $S_1+S_4$ and $S_4+S_6$;
and it is also close to that of $S_6$, which has smallest $x_i$ among $S_1, S_4$ and $S_6$.
The maximum backordering quantity $B^*$ is less than those in the cases of $S_1+S_4$ and $S_4+S_6$;
and it is also close to that of $S_6$, which has smallest $x_i$ among $S_1, S_4$ and $S_6$.
We remark that similar patterns can also be observed in other cases.

\begin{table}
\centering
\scriptsize
\begin{tabular}{lllllllll}
\hline
&$S_1+S_4$&$S_2+S_4$&$S_3+S_4$&$S_1+S_5$&$S_2+S_5$&$S_3+S_5$\\
\hline
$y^*$	&	1630.93	&	1644.44	&	1670.86	&	1649.22	&	1662.63	&	1688.81	\\
$B^*$	&	383.06	&	378.01	&	367.26	&	379.11	&	373.97	&	363.01	\\
$ETPU(y^*,B^*)$	&	1165658.46	&	1160592.66	&	1149976.27	&	1186633.71	&	1181888.2	&	1171942.72	\\
\hline
\\
\hline
&$S_1+S_6$&$S_2+S_6$&$S_3+S_6$&$S_1+S_7$&$S_2+S_7$&$S_3+S_7$\\
\hline
$y^*$	&	1538.28	&	1550.54	&	1574.43	&	1546.29	&	1557.99	&	1580.69	\\
$B^*$	&	207.19	&	201.04	&	187.85	&	192.23	&	185.89	&	172.20	\\
$ETPU(y^*,B^*)$	&	1160290.47	&	1155925.26	&	1146776.05	&	1186440.53	&	1181934.97	&	1172491.58	\\
\hline

\end{tabular}
\caption{Numerical results for different combinations of the screening processes ($n=2$).}
\end{table}

\begin{table}
\centering
\scriptsize
\begin{tabular}{lllllllll}
\hline
&$S_1+S_4+S_6$&$S_2+S_4+S_6$&$S_3+S_4+S_6$&$S_1+S_5+S_6$&$S_2+S_5+S_6$&$S_3+S_5+S_6$\\
\hline
$y^*$	&	1543.69	&	1555.9	&	1579.7	&	1559.95	&	1572.02	&	1595.49	\\
$B^*$	&	205.36	&	199.17	&	185.89	&	199.69	&	193.37	&	179.78	\\
$ETPU(y^*,B^*)$	&	1107457.62	&	1102283.46	&	1091439.66	&	1128854.96	&	1124007.4	&	1113847.66	\\
\hline
\\
\hline
&$S_1+S_4+S_7$&$S_2+S_4+S_7$&$S_3+S_4+S_7$&$S_1+S_5+S_7$&$S_2+S_5+S_7$&$S_3+S_5+S_7$\\
\hline
$y^*$	&	1551.53	&	1563.17	&	1585.73	&	1567.27	&	1578.71	&	1600.83	\\
$B^*$	&	190.33	&	183.94	&	170.15	&	184.43	&	177.89	&	163.74	\\
$ETPU(y^*,B^*)$	&	1131938.47	&	1126598.34	&	1115406.7	&	1154009.02	&	1149005.6	&	1138519.14	\\
\hline

\end{tabular}
\caption{Numerical results for different combinations of the screening processes ($n=3$).}
\end{table}

\section{Concluding Remarks}
In this paper, an inventory model for items with imperfect quality are developed. 
The items are inspected through multiple screening processes and shortage are backordered.
The optimal order size and maximum backordering quantity are obtained by two approaches.
We then provide an analysis of the model. 
Numerical examples are also provided to demonstrate the use of the proposed model.
The current study can be extended in several ways.
In this study, all screening processes are of 100\% accurancy.
The screening rates are also assumed to be greater than the demand rate.
These assumptions may be relaxed in future research.
Another direction is to consider stochastic demand rates.

\newpage

\begin{figure}
\centering 
\begin{tikzpicture} [domain=0:10,xscale=0.8,yscale=0.7] 
\draw[<->] (-0.5,0) -- (-0.5,5);
\node[left] at (-0.5,2.5) {$y$};
\draw (0,0) -- (0,5);
\draw (-0.2,0) -- (0.5,0);
\draw[->] (0.3,0) -- (0.3,1);
\node[right] at (0.3,0.5) {$x_1$};
\node[below] at (0.1,0) {$S_1$};

\draw (2,0) -- (2,0.4) (2,0.6) -- (2,1.4)  (2,1.6) -- (2,2.4)  (2,2.6) -- (2,3.4)  (2,3.6) -- (2,4.4)  (2,4.6) -- (2,5);
\draw (1.8,0) -- (2.5,0);
\draw[->] (2.3,0) -- (2.3,0.8);
\node[right] at (2.3,0.4) {$x_2$};
\node[below] at (2.1,0) {$S_2$};

\draw (4,0) --(4,0.4) (4, 0.6) -- (4,1) (4,1.1) -- (4,1.4) 
(4,1.6) -- (4,2)  (4,2.1) -- (4,2.4)  (4,2.6) -- (4,3)  (4,3.1) -- (4,3.4)
(4,3.6) -- (4,4)  (4,4.1) -- (4,4.4)  (4,4.6) -- (4,5);
\draw (3.8,0) -- (4.5,0);
\draw[->] (4.3,0) -- (4.3,0.6);
\node[right] at (4.3,0.3) {$x_3$};
\node[below] at (4.1,0) {$S_3$};

\node at (5.5,2.5) {$\cdots$};

\draw[dashed] (7,0) -- (7,5);
\draw (6.8,0) -- (7.5,0);
\draw[->] (7.3,0) -- (7.3,0.2);
\node[right] at (7.3,0.1) {$x_n$};
\node[below] at (7.1,0) {$S_n$};

\end{tikzpicture} 
\caption{The screening processes.}
\end{figure}
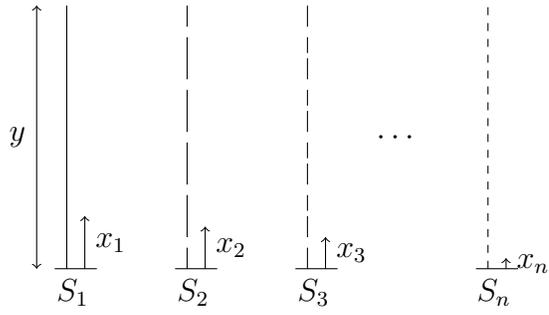 

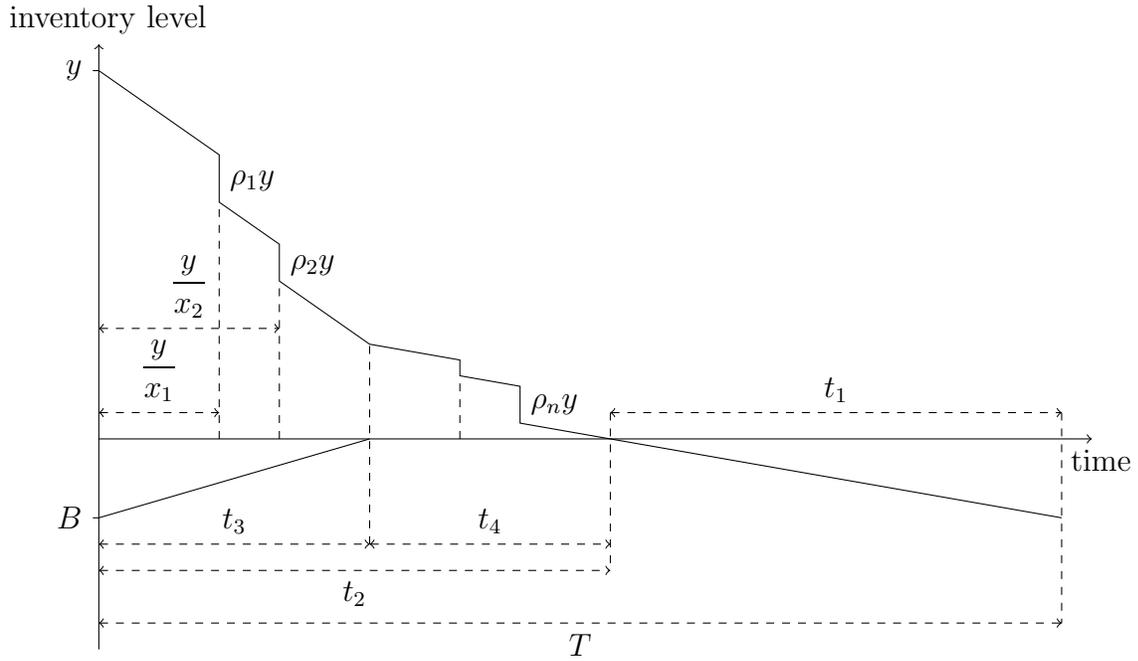
\begin{figure} 
\centering 
\begin{tikzpicture} [domain=0:10,xscale=0.8,yscale=0.7] 
\draw[<->] (0,7.5) node [above] {\; inventory level} -- (0,0) -- (0,0) -- (16.5,0) node [below] {\; time};
\draw (0,0) -- (0,-4);

\draw (0, -1.5) -- (4.5,0);
\draw (0,7) -- (2,5.4) -- (2,4.5) -- (3,3.7) -- (3,3) -- (4.5,1.8) -- (6,1.5) -- (6, 1.2) -- (7, 1) -- (7, 0.3) -- (16,-1.5);
\draw (0,7) -- (-0.1,7) node [left] {$y$};
\draw (0,-1.5) -- (-0.1,-1.5) node [left] {$B$};

\draw[dashed,<->] (0,0.5) -- (2,0.5);
\node[above] at (1,0.5) {$\dfrac{y}{x_1}$};
\draw[dashed,<->] (0,2.1) -- (3,2.1);
\node[above] at (1.5,2.1) {$\dfrac{y}{x_2}$};
\draw[dashed,<->] (0,-2) -- (4.5,-2);
\node[above] at (2.25,-2) {$t_3$};
\draw[dashed] (8.5,0.5) -- (8.5,-2.5);
\draw[dashed,<->] (4.5,-2) -- (8.5,-2);
\node[above] at (6.5,-2) {$t_4$};
\draw[dashed,<->] (0,-2.5) -- (8.5,-2.5);
\node[below] at (4.25,-2.5) {$t_2$};
\draw[dashed,<->] (8.5,0.5) -- (16,0.5);
\node[above] at (12.25,0.5) {$t_1$};
\draw[dashed,<->] (0,-3.5) -- (16,-3.5);
\node[below] at (8,-3.5) {$T$};

\draw[dashed] (2,0) -- (2,4.5);
\draw[dashed] (3,0) -- (3,3);
\draw[dashed] (4.5,-2) -- (4.5,1.8);
\draw[dashed] (6,0) -- (6,1.2);
\draw[dashed] (16,0.5) -- (16,-3.5);

\node[right] at (2,4.9) {$\rho_1 y$};
\node[right] at (3,3.3) {$\rho_2 y$};
\node[right] at (7,0.65) {$\rho_n y$};

\end{tikzpicture} 
\caption{The inventory level in a replenishment cycle.}
\end{figure}

\begin{figure} 
\centering 
\begin{tikzpicture} [domain=0:10,xscale=0.8,yscale=0.7] 
\draw [<->] (0,7.5) node [above] {\; inventory level} -- (0,0) -- (0,0) -- (16.5,0) node [right] {\; time};

\draw (2,0) -- (2,0.9) -- (3,0.9) -- (3,1.6) -- (6,1.6) -- (6,1.9) -- (7,1.9) -- (7,2.6) -- (16,2.6);

\node[left] at (2,0.45) {$\rho_1 y$};
\node[left] at (3,1.3) {$\rho_2 y$};
\node[left] at (7,2.3) {$\rho_n y$};
\node[below] at (2,0) {$\dfrac{y}{x_1}$};
\draw[dashed] (3,0) -- (3,0.9);
\node[below] at (3,0) {$\dfrac{y}{x_2}$};
\node[below] at (7,0) {$\dfrac{y}{x_n}$};
\draw[dashed] (7,0) -- (7,1.9);
\node[below] at (16,0) {$T$};
\draw[dashed] (16,0) -- (16,2.6);


\end{tikzpicture} 
\caption{The inventory level of defective items.}
\end{figure}
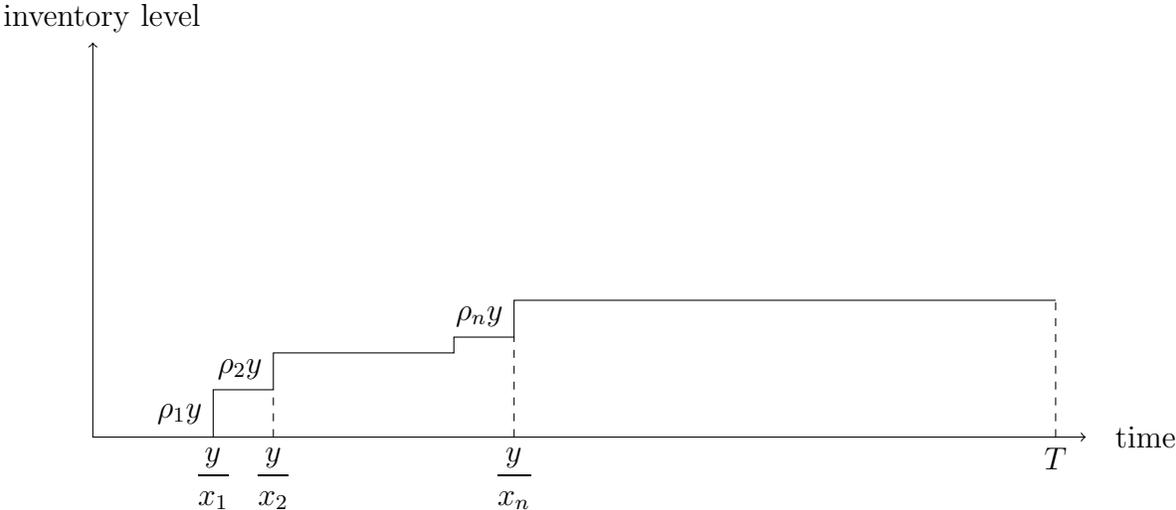 
\end{document}